\documentclass[a4paper,12pt]{article}

\usepackage{amssymb}
\usepackage{textcomp}
\usepackage{tipa}

\usepackage[T2A]{fontenc}
\usepackage[cp1251]{inputenc}
\usepackage[tbtags]{amsmath}

\usepackage[dvips]{graphicx}
\usepackage[english,russian]{babel}

\usepackage{wasysym, amsthm
}

\theoremstyle{definition}

\newtheorem{definitionhead}{Определение}

\begin{document}
\centerline{\LARGE Алгоритмические методы и комбинаторика слов}
\centerline{\LARGE в теории колец}
\centerline{\large (Черновик)}

\medskip
\centerline{\Large М.И. Харитонов
\footnote{М. И. Харитонов. E-mail: mikhailkharitonov@yandex.ru}}

\medskip
\begin{itemize}

	\item Приведённые ниже факты отдельно доказывались для $\epsilon$- и $Lie$- алгебр, но доказательства и формулировки для различных типов алгебр похожи, поэтому ниже приведена попытка объединения формулировок.

А.И. Ширшов ввёл и доказал следующие понятия и теоремы:

{\bf Определение} Слова длины 1 назовем {\em $\Omega$-правильными}
($\Omega = K, AK, Lie$) словами и произвольно упорядочим. Считая, что $\Omega$-правильные слова, длина которых меньше $n$, $n >1$, уже определены и упорядочены каким-то способом так, что слова меньшей длины предшествуют
словам большей длины, назовем слово $w$ длины $n$ {\em $\Omega$-правильным}, если
\begin{enumerate}
 \item $w =uv$, где $u$, $v$ — $\Omega$-правильные слова;
 \item $u\geqslant v$ при $\Omega = K$ и $и>v$ при $\Omega = AK, Lie$;
 \item (Только для $\Omega=Lie$) если $u=u_1u_2$, то $u\leqslant v$.
\end{enumerate}

{\bf Теорема 1.} Правильные слова образуют базис свободной $\Omega$-алгебры.

{\bf Теорема 2.} Всякая подалгебра $\Psi$ свободной $\Omega$-алгебры $\Phi$ свободна. 

{\bf Проблема тождества для $\Omega$-алгебр.} Существует ли алгоритм, который для произвольного конечного множества $S$ и произвольного элемента $a$ из $\Omega$-алгебры позволяет выяснить, принадлежит ли $a$ идеалу $\langle S\rangle$.

{\bf Теорема о тождестве 1.} Пусть $S$ -- некоторое фиксированное множество элементов свободной $\Omega$-алгебры $E$. Тогда существует алгоритм, позволяющий за конечное число шагов определить, принадлежит ли произвольный элемент $t\in E$ идеалу $\langle S\rangle$.

{\bf Следствие.} Существует алгоритм, решающий проблему тождества для алгебр Ли с одним определяющим соотношением.

{\bf Теорема о тождестве 2.} Существует алгоритм, решающий проблему тождества для алгебр Ли с однородными множествами определяющих соотношений.

{\bf Теорема о свободе.} Пусть $E_0$ -- $\Omega$-алгебра с множеством порождающих $R$ и одним опредлеяющим соотношением $s=0$, в левую часть которого входит образующий $a_\alpha$. Тогда подалгебра $E'_0$, порождённая в алгебре $E_0$ множеством $R\setminus a_\alpha$, свободна.

Литература: \cite{Sh53}, \cite{Sh54}, \cite{Sh58}, \cite{Sh62(1)}, \cite{Sh62(2)}, \cite{BBL97}. 

	\item Проблема Куроша-Левицкого для конечно порождённых 
\begin{itemize}
	\item ниль-алгебр конечного ниль-индекса
	\item алгебр конечного индекса
\end{itemize}
{\bf Теорема Ширшова о высоте.} Множество всех не $n$-разбиваемых слов в конечно порождённой алгебре с допустимым полиномиальным тождеством имеет ограниченную
 высоту $H$ над множеством слов степени не выше $n-1$.

Литература: \cite{Kur41}, \cite{BBL97}, \cite{SGS78}.

	\item {\bf Определение.} Ассоциативное слово называется {\em правильным}, если оно лексикографически больше любого своего циклического сдвига.

Неассоциативное слово называется {\em правильным}, если оно правильное в ассоциативном смысле и 
\begin{itemize}
	\item если $[u] = [[v][w]]$, то $v$ и $w$ -- правильные слова,
	\item если $[u] = [[v_1] [v_2]] [w]$, то $v_2 \leqslant w$.
\end{itemize}

{\bf Теорема Ширшова}. В правильном в ассоциативном смысле слове существует единственный способ расставить Лиевы скобки так, чтобы полученное слово было правильным в неассоциативном смысле. 

Правильные слова образуют базис свободной алгебры Ли.

Литература: \cite{BBL97}, \cite{SGS78}, \cite{Sh58}, \cite{Ufn89}. 

	\item {\bf Определение.} Назовём слово $а$ полуправильным, если
любой его конец либо лексикографически меньше $a$, либо является
началом $a$.

{\bf Теорема.} Любое бесконечное слово над конечным алфавитом содержит подслово $fgf$, где $f$ -- полуправильное, а $g$ -- правильное (возможно, пустое) слово.

Литература: \cite{BBL97},  \cite{Ufn89}. 

	\item {\bf Теорема Ван дер Вардена.} Пусть $n$ и $k$ -- натуральные числа, последовательность натуральных чисел разбита на $k$ множеств. Тогда найдётся число $f(n, k)$ такое, что среди первых $f(n, k)$ натуральных чисел найдётся арифметическая прогрессия длины $n$ из одного множества. 

Многомерное обобщение для фигур и гомотетии с положительным коэффициентом.

Литература: \cite{Bug06}, \cite{Hin79}. 

	\item {\bf Определение.} Группа удовлетворяет условию $C'(\lambda)$, когда общая часть любых двух порождающих соотношений меньше, чем $\lambda$, умноженное на длину любого из них.

{\bf Лемма Гриндлингера.} В карте, удовлетворяющей условию $C'({1\over 6})$, найдётся клетка, большая часть границы которой лежит на границе карты.

Алгебраическая формулировка с группами и соотношениями.

Алгоритм Дена-Гриндлингера определения тривиальности группового слова в группе с конечным числом соотношений.

Литература: \cite{Sap14}, \cite{Kl09}. 

	\item {\bf Теорема Регева.} Если алгебры $A$ и $B$ удовлетворяют полиномиальному тождеству, то алгебра  $A\otimes_F B$ также удовлетворяет полиномиальному тождеству. 

Литература: \cite{Reg71}, \cite{Lat72}. 

	\item {\bf Diamond-lemma.} Пусть $M$ -- ЧУМ, в котором любая убывающая цепь -- конечна. 

{\bf Определение.} {\em Отношение Чёрча-Россера:} $x\leftrightsquigarrow y$, если у $x$ и $y$ есть общий потомок.

Представим $M$ в виде графа Ньюмана с множеством рёбер $R$. Тройка $(M, \leqslant, R)$ называется {\em схемой симплификации}. Следующие условия эквивалентны:
\begin{enumerate}
\item $M$ -- обладает свойством каноничности (т.е. у каждого $m\in M$ нормальная форма единственна).
\item Отношение Чёрча-Россера -- транзитивно.
\item Выполняется условие локального слияния (``у любых двух братье есть общий потомок'').
\item В любой компоненте связности лежит ровно один минимальный элемент.
\item ($x\backsim y$, т.е. между $x$ и $y$ есть неориентированный путь) $\Longleftrightarrow (x\leftrightsquigarrow y)$.
\end{enumerate}

\medskip
{\bf Определение.} Введём на мономах $X^*$ линейный порядок $<$ такой, что для любого монома $z\in X^*$ имеет место $x<y\Rightarrow xz<yz$. {\em Базис Грёбнера-Ширшова} некоторого идеала $I\vartriangleleft  k\langle X\rangle$ -- это конечное множество полиномов $G$, порождающее идеал $I$, причём старший моном $\bar{h}$ любого полинома $h\in I$ делится на некоторый старший моном $\bar{g}$ полинома из базиса Грёбнера-Ширшова.

Элемент $h$ обладает {\em $H$-представлением} относительно системы порождающих $G$, если в представлении $h=\sum \alpha_i u_ig_iv_i$ любой моном $u_i\bar{g_i}v_i$ не больше, чем $\bar{h}$.

Определим для полинома $f\in k\langle X\rangle$ его $supp(f)$ -- упорядоченное множество составляющих его мономов. Тогда лексикографический порядок на суппортах полиномов индуцирует частичный порядок $\leqslant_{supp}$ на полиномах $k\langle X\rangle$.

{\bf Теорема.} Следующие условия эквивалентны.

\begin{enumerate}
\item $G$ -- базис Грёбнера-Ширшова $I$.
\item Любой элемент $I$ редуцируется относительно $G$ к нулю.
\item Любой $h\in I$ обладает $H$-представлением относительно $G$.
\item Схема  $(k\langle X\rangle, \leqslant_{supp}, R_G)$ обладает свойством каноничности.
\end{enumerate}

Литература: \cite{Lat12}, \cite{Ber78}, \cite{BMM95}.

	\item {\bf Теорема (Туэ -- 1).} Пусть $X = \{a, b\}$, подстановка $\phi$ задана соотношениями $\phi(a)=ab$, $\phi(b)=ba$. Тогда если слово $w\in X^*$ -- бескубное, то и $\phi(w)$ -- бескубное.

{\bf Теорема (Туэ -- 2).} Пусть $X = \{a, b, c\}$, подстановка $\phi$ задана соотношениями $\phi(a)=abcab$, $\phi(b)=acabcb$, $\phi(c)=acbcacb$. Тогда если слово $w\in X^*$ -- бесквадратное, то и $\phi(w)$ -- бесквадратное.

{\bf Теорема (Туэ -- 3).} Пусть $M$ и $N$ -- алфавиты, для подстановки $\phi:M\rightarrow N^*$ и выполнены следующие условия:
\begin{enumerate}
\item если длина $w$ не больше 3, то $\phi(w)$ -- бесквадратное;
\item если $a$, $b$ -- буквы алфавита $M$, а $\phi(a)$ -- подслово $\phi(b)$, то $a=b$. 
\end{enumerate}
Тогда если слово $w\in M^*$ -- бесквадратное, то и $\phi(w)$ -- бесквадратное.

{\bf Теорема (Крошмор).} Пусть $\phi$ -- подстановка, $M$ -- наибольший размер блока, $m$ -- минимальный размер блока, $k=max\{3, 1+[(M-3)/m]\}$. Тогда подстановка $\phi$ -- бесквадратная в том и только в том случае, когда для любого бесквадратного слова $w$ длины $\leqslant k$ слово $\phi(w)$ будет бесквадратным.

Литература: \cite{Sap14}.

	\item  {\bf Определение.} Алгебра $А$ называется {\em мономиальной}, если в ней есть базис определяющих
соотношений вида $c = 0$, где $c$ — слово от образующих алгебры.

{\em Конечным автоматом ({\bf КА})} с алфавитом $X$
входных символов называется ориентированный граф, в котором выделено два (возможно
пересекающиеся) множества вершин, называемых {\em начальными}
и {\em финальными (конечными)} и каждое ребро помечено буквой из конечного
алфавита $X$. 
Язык $L$ называется
{\em регулярным} или {\em автоматным}, если существует конечный
автомат, допускающий слова из множества $L$ и только их.

Автомат называется {\em детерминированным}, если 
\begin{enumerate}
\item начальная вершина ровно одна;
\item из любой его вершины не может выходить более одного ребра, помеченного
одной и той же буквой;
\item нет ребер, помеченных пустой цепочкой.
\end{enumerate}

{\bf Предложение.} для всякого недетерминированного
{\bf КА} существует детерминированный {\bf КА},
допускающий то же самое множество слов.

{\bf Определение.} Алгебра $A$ называется {\em автоматной},
если множество ее ненулевых слов от образующих
А является регулярным языком.

{\bf Предложение.} Конечно определенная мономиальная
алгебра является автоматной.

{\bf Определение.} {\em Функция роста $V_A(n)$ алгебры} $A$ — это размерность пространства,
порожденного словами длины не выше $n$.

Если следующий предел существует,
то его значение называется {\em размерностью Гельфанда—Кириллова}
алгебры $A$ и обозначается $GK(A)$:
$$GK(A)=\lim\limits_{n\rightarrow\infty}{{\ln(V_A(n))\over \ln(n)}}.$$

Пусть $\Gamma(A)$ — минимальный детерминированный граф автоматной
алгебры $A$. Назовем вершину графа {\em циклической}, если
существует путь, начинающийся и заканчивающийся в этой
вершине. Назовем вершину {\em дважды циклической}, если существуют
два различных пути, начинающихся и заканчивающихся
в этой вершине и не проходящих ни через одну другую вершину
дважды.

Пусть граф $\Gamma$ не имеет дважды циклических вершин. Назовем
{\em цепью} подграф графа $\Gamma$, состоящий из последовательности
ребер, в которой конец предыдущего ребра является началом
следующего, и никакая вершина не встречается дважды. Назовем
{\em простым графом} подграф графа $\Gamma$, состоящий из конечного
числа циклов, занумерованных числами $1, 2,\dots, d$, причем пары
соседних циклов с номерами $i, i + 1$ соединены ровно одной
цепью, направленной от $i$-ro к $(i+1)$-му циклу. В первый цикл может
входить одна цепь, и из последнего также может выходить
одна цепь. 

{\bf Теорема (Уфнаровский).} Пусть $A$ — автоматная
алгебра, $\Gamma(A)$ — ее минимальный детерминированный граф.
\begin{enumerate}
\item Если $\Gamma(A)$ имеет вершину, принадлежащую двум различным
циклам, то $A$ имеет экспоненциальную функцию роста.
\item Если $\Gamma(A)$ не имеет дважды циклических вершин,
то $A$ имеет степенную функцию роста. Степень роста
(размерность Гельфанда—Кириллова) равна количеству циклов
в максимальном простом подграфе, содержащемся в $\Gamma(A)$.
\end{enumerate}

{\bf Теорема.} Пусть граф автоматной мономиальной алгебры
$A$ не имеет вершин, принадлежащих двум циклам. Тогда
$A$ вкладывается в алгебру матриц над полем.

{\bf Следствие.} Пусть $A$ — автоматная мономиальная
алгебра, $\Gamma(A)$ — ее минимальный детерминированный граф. Тогда следующие уловия эквивалентны:
\begin{enumerate}
\item $\Gamma(A)$ не имеет дважды циклических вершин;
\item алгебра $A$ имеет степенной рост;
\item алгебра $A$ имеет не экспоненциальный рост;
\item алгебра $A$ представима матрицами над полем;
\item в $A$ выполняется полиномиальное тождество.
\end{enumerate}
Литература: \cite{BBL97}, \cite{Sap14}, \cite{Ufn89}, \cite{wMA}, \cite{wDS}.
\end{itemize}

\end{document}